\documentstyle[12pt]{article}

\def\A{{\cal A}}
\def\B{{\cal B}}
\def\Cc{{\cal C}} 
\def\dim{{\rm dim}}
\def\E{{\cal E}} 
\def\F{{\cal F}} 
\def\G{{\cal G}} 
\def\GG{{\bf G}}
\def\gg{{\bf g}}
\def\lra{{\longrightarrow}}
\def\O{{\cal O}} 
\def\pp{{\bf p}}
\def\V{{\cal V}}

\newtheorem{theorem}{Theorem}[section]
\newtheorem{proposition}[theorem]{Proposition}
\newtheorem{lemma}[theorem]{Lemma}

\def\num{\refstepcounter{theorem}\paragraph{ (\thetheorem)}}
\def\rem{\refstepcounter{theorem}\paragraph{Remark \thetheorem}}
\def\proof{\paragraph{Proof}}

\makeatletter \def\l@section{\@dottedtocline{1}{0em}{1.2em}} \makeatother

\begin{document}

\title{Boundedness of semistable principal bundles on a curve, 
with classical semisimple structure groups}
\author{Nitin Nitsure}
\date{18 February 1998}
\maketitle

\setcounter{section}{-1}
\section{Introduction}
Let $k$ be an algebraically closed field of characteristic $\ne 2$.
Let $G$ be a semisimple, simply connected algebraic group over $k$
which is of classical type, that is, $G$ is special linear,
special orthogonal, or symplectic. 
Let $X$ be a nonsingular irreducible projective curve over $k$.
We consider \'etale locally trivial 
principal $G$-bundles on $X$ which are semistable.
The main result in this paper is the following:

{\bf Theorem 1 } {\sl 
If $G$ is a semisimple, simply connected algebraic group of 
classical type then semistable principal $G$-bundles
on $X$ form a bounded family.
}

In other words, we prove that there exists a family of principal 
$G$-bundles on $X$ parametrized by a scheme $T$ of finite type 
over $k$ such that every semistable $G$-bundle on $X$ is isomorphic
to the restriction of this family to $\{t\}\times_k X$ for 
some closed point $t$ of $T$. 

In characteristic zero, the theorem follows from the 
Narasimhan-Seshadri theorem, and $G$
can be any semisimple group, not necessarily of classical type
(see Ramanathan [R]).
The proof given below is characteristic free, 
but needs restriction to classical type.

The theorem is already known when
$G$ is a special linear group $SL_n$, 
as in that case a semistable principal
$G$-bundle on $X$ `is the same as' a semistable vector bundle
on $X$ whose rank is $n$ and determinant is trivial. 
We therefore have to prove the theorem only when $G$ is special 
orthogonal group $SO_n$ or symplectic group $Sp_{2n}$. 
This is done as follows.

A principal $G$-bundle in these cases is `the same as' a pair
$(V,b)$ where $V$ is a vector bundle on $X$ with trivial determinant
and $b:V\otimes_{\O_X}V\to \O_X$ is a non-degenerate bilinear form,
which is either symmetric or skew-symmetric. 
Let $F\subset V$ be a totally
isotropic vector subbundle of the maximal possible rank $r$
(isotropic means the restriction 
$b:F\otimes_{\O_X}F\to \O_X$ is identically zero, 
and the maximal possible rank is the rank of $G$). 
The crucial step is to show that 
the maximal possible degree of such subbundles
is bounded below, where the bound depends only on $X$ and $G$. 
This we do by a generalization of the 
Mukai-Sakai theorem (see [M-S]) 
for vector bundles to the present situation.
On the other hand, if the principal $G$-bundle corresponding
to $(V,b)$ is semistable, then an elementary calculation shows that 
the degree of any isotropic subbundle is less than or equal to zero.

It follows that the set of all isotropic subbundles $F$ 
of maximal rank and maximal degrees of all possible 
$(V,b)$ corresponding to semistable $G$-bundles (where $G$ is fixed)
constitutes a bounded family of vector bundles. From this we
can deduce that the semistable $G$-bundles
form a bounded family.

The arrangement of this paper is as follows. In section 1 we recall
the basic definitions, and prove that a pair $(V,b)$ is semistable
(if and) only if every isotropic subbundle
of $(V,b)$ has degree $\le 0$.
In section 2, we prove an analogue of the Mukai-Sakai theorem,
showing that the maximal possible degree of isotropic subbundles
is not too small. In section 3 we complete the proof of theorem 1. 

{\bf Acknowledgements } I thank M.\,S.\,Narasimhan for bringing the 
boundedness problem to my (and general) attention. 
I thank M.\,S.\,Raghunathan for sharing his
view that a solution lay in the direction of parabolic reductions; this
view is vindicated by the proof presented here. 
I thank T.\,N.\,Venkataramana for answering my questions on the subject
of semisimple groups.

\section{Semistability and isotropic subbundles}

Let $E$ be a principal $G$-bundle, and $P\subset G$ a maximal
parabolic subgroup. Let $\sigma:X\to E/P$ be a section, which is
the same as a reduction of the structure group to $P$.
Let $T_{\pi}$ be the relative tangent bundle to the projection
$\pi :E/P \to X$. Recall (see Ramanathan [R]) that the bundle
$E$ is said to be {\bf semistable} if for every maximal parabolic
$P\subset G$ and every reduction $\sigma:X\to E/P$, the 
following inequality holds:
\num $~~~~~~~~~~~~~~~~~~~~~~~~~~~
\deg (\sigma^*(T_{\pi})) \ge 0 $

\medskip

If $G=SL_n$, then the principal bundle $E$ corresponds to
vector bundle $V$ with trivial determinant. A reduction to a 
maximal parabolic corresponds to giving a vector subbundle
$F\subset V$, and the above inequality is then equivalent
to the usual definition of semistability for vector bundles. 
Next we consider the case when $G$ is a special orthogonal 
or symplectic group. Let $W$ be the vector space $k^{2n}$
(or $k^{2n+1}$), with standard basis $(e_i)$. We put
$u_i=e_i$ and $v_i=e_{i+n}$ for $1\le i\le n$ 
(and $w=e_{2n+1}$). Consider the 
bilinear form $b$ defined in one of the following ways.

For the group $SO_{2n+1}$, we put 
\num $~~~~~~~~~~~~~~~
b(u_i,v_i)=b(v_i,u_i)=b(w,w)=1$ for $1\le i\le n$,

\medskip

while all other pairings of basis vectors are zero.

For the group $Sp_{2n}$, we put 
\num $~~~~~~~~~~~~~~~
b(u_i,v_i)= - b(v_i,u_i)=1$ for $1\le i\le n$,

\medskip

while all other pairings of basis vectors are zero.

For the group $SO_{2n}$, we put 
\num $~~~~~~~~~~~~~~~
b(u_i,v_i)=b(v_i,u_i)=1$ for $1\le i\le n$,
 
\medskip

while all other pairings of basis vectors are zero.

Let $G\subset SL(W)$ be the subgroup preserving $b$. 
Let $P_r\subset G$ be the subgroup which carries the $r$-dimensional
linear subspace $W_r = <u_1,\ldots,u_r>$ into itself.
Then $P_r$ is a maximal parabolic subgroup of $G$. Let $T\subset G$
be the standard maximal torus, consisting of matrices of the form
$$diag(a_1,\ldots, a_n,a_1^{-1},\ldots, a_n^{-1})$$ 
when $W$ is even dimensional (that is, $G$ is $SO_{2n}$ or $Sp_{2n}$), 
and of the form
$$diag(a_1,\ldots, a_n,a_1^{-1},\ldots, a_n^{-1}, 1)$$
when $W$ is odd dimensional (that is, $G$ is $SO_{2n+1}$).
Let for $1\le i \le n$, the multiplicative character
$$\lambda_i :T\to \GG_m$$
be defined to take the value $a_i$ on the above diagonal matrices.

Let $\gg$ and $\pp_r$ be the Lie algebras of 
$G$ and $P_r$ respectively.
Then $T$ acts by adjoint action on the quotient $\gg/\pp_r$, and
$T$ acts on $W_r$ by restriction of the defining representation 
of $G$ on $W$. 

The following lemma can be proved by an elementary calculation,
which we omit.
\begin{lemma}\label{C(n,r)} 
(i) The torus $T$ acts on $\det(W_r)$ by the character
$\alpha_r =\lambda_1 + \ldots +\lambda_r$. 

(ii) The torus $T$ acts on $\det(\gg/\pp_r)$ by the  
character $\chi_r$, given case by case as follows.
$$ \chi_r =  \left\{ 
\begin{array}{ll}
-(2n-r)\alpha_r & \mbox{if $G=SO_{2n+1}$} \\
 -(2n-r+1)\alpha_r & \mbox{if $G=Sp_{2n}$} \\
-(2n-r-1)\alpha_r & \mbox{if $G=SO_{2n}$ }
\end{array}\right. $$

(iii) In particular, in each case $\chi_r$ is a negative 
multiple of the character $\alpha_r$ on $\det(W_r)$,
which we write as

\centerline{$\chi_r = - C(n,r)\alpha_r$}

where $C(n,r)>0$ is an integer given as above. 
\end{lemma}

Now let $E$ be a principal $G$-bundle, where $G$ is as above.
Let $(V,b)$ be the associated vector bundle together with
a bilinear form $b$. Then for any isotropic subbundle
$F_r\subset V$ of rank $r$, we get a reduction of structure
group to (upto isomorphism) the parabolic subgroup $P_r$.
If the transition functions are contained in the torus $T$,
then it follows from the above lemma that the associated
$\gg/\pp_r$-bundle $\F_r$ has degree equal to the strictly negative
multiple 
\num\label{multiple}
$~~~~~~~~~~~~~~~~~~\deg(\F_r) = -C(n,r)\deg(F_r)$ 

\medskip

of the degree of $F_r$. The above relation (\ref{multiple})
holds even if the structure group is not reduced to $T$, 
as we can have a flat deformation which reduces the 
structure group to $T$, and such a deformation does not 
affect the degrees.

Note that the bundle $\F_r$ is canonically isomorphis to
the pullback $\sigma^*(T_{\pi})$ of the relative tangent bundle
of $\pi :E/P \to X$ by the parabolic reduction $\sigma$.
Hence if $E$ is semistable, we must have $\deg(\F_r)\ge 0$.
Hence the equation (\ref{multiple}) now completes
the proof of the `only if' part of the following proposition.

\begin{proposition}\label{negative}
The principal bundle corresponding to a pair $(V,b)$ 
is semistable if and only if $\deg(F)\le 0$ 
for any isotropic subbundle $F$ of $(V,b)$.
\end{proposition}

For the `if' part, just observe that any maximal parabolic
subgroup of $G$ is conjugate to one of the subgroups $P_r$
that we have explicitly described.

\rem In particular, if $E$ is semistable, then any 
subbundle of an isotropic subbundle $F$ of $(V,b)$ 
has non-positive degree, so the Harder-Narasimhan type 
of $F$ is bounded above.

\section{Isotropic subbundles of maximal degrees}

In this section, we generalize the result of Mukai-Sakai [M-S]
to isotropic subbundles of $(V,b)$. The method is a
straightforward generalization of [M-S].
In this section, we do not assume that $(V,b)$ is stable.

We continue to use the same notation as the previous section.
Let $(V,b)$ correspond to a principal $G$-bundle as above.
Let $F_r\subset V$ be an isotropic subbundle of rank $r$
such that $F_r$ has the maximum possible degree (say $d$)
amongst isotropic subbundles of $V$ having the same rank $r$.
Note that by Riemann-Roch theorem for the curve $X$,
the Hilbert polynomial of $V/F_r$ is
$$h(t) = -d + (n-r)(1-g_X) + (n-r)t$$ 
where $g_X$ denotes the genus of $X$.
Let $Q$ be the quot scheme parametrizing all equivalence
classes of quotients 
$$q: V \to \E$$
where $\E$ is a coherent sheaf with Hilbert polynomial 
$h(t)$. Then on $X\times _k Q$ we have a universal
quotient 
$$q: p_1^*V \lra \G$$
Then $Q$ is projective over $k$. 
By pulling back $b$ under $X\times _kQ \lra X$, we
get a bilinear form $p_1^*(b)$ on $p_1^*V$.
Let $Q^{iso}\subset Q$ be the closed subscheme where the 
kernel of $q: p_1^*V \lra \E$ is isotropic. Let 
$Q_0^{iso}$ be an irreducible component of $Q^{iso}$
which contains the quotient $V\lra V/F_r$ that we started with.

\begin{lemma} If $F_r$ has maximal degree amongst all rank $r$ 
isotropic subbundles in $V$, then the restriction
$$\G | X\times _k Q_0^{iso}$$
of the universal quotient $\G$ to any irreducible component
$Q_0^{iso}$ containing $V\lra V/F_r$ is a locally free sheaf
on $X\times _k Q_0^{iso}$.
\end{lemma}

\proof Suppose not, then there exists a closed point 
$t\in Q_0^{iso}$ such that the restriction 
$\G_t \,=\,\G | X\times _k \{ t\}$
is not locally free. Let $F$ be the kernel of $V\lra \G_t$.
Let $F'$ be the $\O_X$-saturation of $F$. Then $F'$ is
generically equal to $F$, so $F'$ is isotropic and of the 
same rank $r$. But as $\deg(F') > \deg(F)$, this contradicts 
the maximality of the degree of $F_r$ as $\deg(F)=\deg(F_r)$.

The following proposition generalizes the theorem of Mukai-Sakai
to isotropic subbundles.

\begin{proposition}\label{nottoosmall}
Let $F_r$ be an isotropic subbundle of rank $r$ in $(V,b)$
such that the degree of $F_r$ is maximal amongst such subbundles.
Then the following inequality holds
$$\deg(F_r) \ge { -g_X \,\cdot\,\dim(G/P_r)  \over C(n,r)}$$
where $C(n,r)$ is the positive integer given by lemma \ref{C(n,r)}
and $G/P_r$ is the quotient of $G$ by the maximal parabolic $P_r$,
and $g_X$ is the genus of $X$. 
\end{proposition}

\proof As $G$ acts on the left on $G/P_r$, the principal $G$-bundle
$E$ has an associated bundle $Y_r \lra X$ with fiber $G/P_r$.
Note that $Y_r$ is the same as the closed subscheme, defined
by the condition of isotropy, of the Grassman bundle of
rank $r$ subspaces of the fibers of $V$.

Let $Q_0^{iso}$ be an irreducible component of $Q^{iso}$ containing
$V\lra V/F_r$. Then by the above lemma, we have a short
exact sequence of vector bundles
$$0 \lra F \to p_1^*V \lra \G \lra 0$$
on $X\times _k Q_0^{iso}$. Hence we get a morphism
$$\varphi : X\times_k Q_0^{iso} \lra  Y_r$$
over $X$, which is the classifying morphism for the isotropic 
subbundle $F \subset p_1^*V$. 

Note that as both sides are projective, the morphism
$\varphi$ is proper. We claim that in fact $\varphi$ is finite. 
For, as in the corresponding argument of Mukai-Sakai, 
otherwise there will exist a complete curve $B\subset Q_0^{iso}$
and a closed point $x\in X$ such that $\{ x \} \times_k B$
lies in a fiber of $\varphi$. Then by rigidity,
the restricted morphism  
$$\varphi : X \times_k B \lra Y_r$$
will factor through the projection $X \times_k B \lra B$,
giving a contradiction, as $\varphi$ is over $X$. 
This shows that $\varphi$ is finite, hence 
\num\label{half}$~~~~~~~~~~~~~~~
\dim (Q_0^{iso}) \le \dim (G/P_r)$

\medskip

Note that each $S$-valued point of $Q_0^{iso}$, where $S$ is a
$k$-scheme, can be regarded as a section $\sigma : X\times_kS
\lra Y_r\times _kS$. For any $k$-scheme $S$ and any section $s$ of 
$Y_r\times_kS\lra X\times_kS$, we
get an isotropic subbundle of the pullback of $V$ to $X\times_kS$.
This gives a morphism $U \lra Q_0^{iso}$
where $U$ is the scheme of sections of $Y_r\to X$, which
is clearly injective on $k$-valued points. 
It follows that the dimension of
$Q_0^{iso}$ is greater than or equal to the maximum of the
dimensions of irreducible components $U_0$ of $U$ which contain 
the section $\sigma_r :X\lra Y_r$ which corresponds to the 
subbundle $F_r\subset V$. 

We now quote the following proposition, due to Mori.
 
\begin{proposition}\label{Mori}
Let $\pi: Y\lra X$ be a projective morphism, and let $U$ be the 
scheme of all sections of $Y\lra X$. Let $T_{\pi}$ be the 
relative tangent bundle to $\pi:Y\lra X$, and let 
$\sigma : X \lra Y$ be a closed point of $U$.
Let $U_0$ be an irreducible component of $U$ of maximal 
dimension which contains $\sigma$. 
Then the following inequality holds.
$$\dim(U_0) \ge \dim H^0(X, \sigma^*(T_{\pi})) -
\dim H^1(X, \sigma^*(T_{\pi}))$$ 
\end{proposition}

\proof This follows from proposition 3 of Mori [M], by
taking $Z$ to be empty in the notation of [M], and observing  
that $U$ is an open subscheme of the scheme 
$\underline{Hom}(X,Y)$.

We now apply this to the present case, by taking $X$
to be our curve $X$, $Y$ to be $Y_r$, and $\sigma$ to be
the section $\sigma_r$ corresponding to the maximal
degree isotropic subbundle $F_r$. Then 
$\sigma^*(T_{\pi})$ equals the associated $\gg/\pp_r$ 
bundle $\F_r$ in the notation of section 1.
By equation \ref{multiple}, 
we have $\deg(\F_r) = -C(n,r) \cdot d$
Hence by Riemann-Roch for the curve $X$ we get
$$\dim H^0(X, \sigma^*(T_{\pi})) -
\dim H^1(X, \sigma^*(T_{\pi})) =
-C(n,r)\cdot d + \dim(G/P_r)(1-g_X)$$
Combining the above equation with proposition \ref{Mori},
we get
$$-C(n,r) \cdot d\, +  \dim(G/P_r)(1-g_X) \le \dim(U_0)$$
Now, $U_0$ embeds in $Q_0^{iso}$, hence combining the above
with equation \ref{half}, we finally have
$$-C(n,r) \cdot d + \dim(G/P_r)(1-g_X) \le \dim(G/P_r)$$
Solving this for $d$ 
completes the proof of the proposition \ref{nottoosmall}.

\section{Proof of theorem 1}

We first recall the following standard fact.

\rem\label{ext} Let $\A$ be a set of isomorphism classses 
of vector bundles on $X$. We say that $\A$ is {\bf bounded}
if there exists 
a finite type scheme $T$ over $k$ and vector bundles
$\E$ on $X\times_k T$ such that given any element $a\in \A$, there
exists a closed point $t\in T$ such that the restriction
$\E|(X \times_k \{ t\})$ represents $a$. Suppose $\A$ and $\B$
are two sets of isomorphism classes of vector
bundles on $X$. Let $\Cc$ be the set of all isomorphism
classes of vector bundles $\V$ which fit in a short exact
sequence
$$0 \lra \E \lra \V \lra \F \lra 0$$
where $[\E]\in \A$ and $[\F]\in \B$. Then it is known that
if $\A$ and $\B$ are bounded sets, then $\Cc$ is also a bounded 
set. 

With this, we are ready to prove theorem 1.

\proof The above remark implies that the set of all 
isomorphism classes of rank $n$ vector bundles $F$ on $X$ which 
occur as a maximal degree isotropic vector subbundle of $(V,b)$
where $(V,b)$ corresponds to a semistable principal
$G$-bundle on $X$ (where $G$ is one of $SO_{2n+1}$, $Sp_{2n}$,
$SO_{2n}$) is bounded, as the degrees of such $F$ are
bounded below by proposition \ref{nottoosmall}, 
and their Harder-Narasimhan
types are bounded above by proposition \ref{negative}.
Note that as $F$ is isotropic of the maximal possible rank $n$, 
the bilinear form $b$ induces an isomorphism $F^* \cong V/F$.
The remark \ref{ext} implies that the set of isomorphism classes
of extensions of the type 
$$0 \lra F \lra \E \lra F^* \lra 0$$
form a bounded set. This shows that the isomorphism classes of the 
underlying vector bundles $V$ of semistable
pairs $(V,b)$ form a bounded set. 

Let $S$ be a finite type scheme over $k$ and
$\V$ a vector bundle on $X\times S$ in which all
such $[V]$ occur. Hence there is a linear scheme
$T$ over $S$ which parametrizes the 
pairs $(V,b)$ where $b$ is a bilinear form on $V$.
Then $T$ has a locally closed subscheme where $b$ is
non-degenerate, and symmetric (or skew symmetric).
This completes the proof that the isomorphism classes 
of semistable $(V,b)$ form a bounded set.

\section*{References} \addcontentsline{toc}{section}{References}

[M] S.\,Mori : Projective manifolds with ample tangent bundles.
Ann. of Math. {\bf 110} (1979) 593-606.  

[M-S] S.\,Mukai and F.\,Sakai : Maximal subbundles of vector bundles
on a curve. \hfill\\
Manuscripta math. {\bf 52} (1985) 251-256.

[R] A.\,Ramanathan : Stable principal bundles on a compact 
Riemann surface. \hfill\\
Math. Annln. {\bf 213} (1975) 129-152.

\bigskip
Address:

School of Mathematics, Tata Institute of Fundamental Research,
Homi Bhabha Road, Mumbai 400 005, India. e-mail:
nitsure@math.tifr.res.in

\end{document}